\documentclass[reqno]{amsart}

\usepackage{graphicx} 
\usepackage{mathdots}
\usepackage{tikz-cd}
\usepackage{bbm}
 
\usepackage{amsmath,amssymb,latexsym, amsfonts, amscd, amsthm} 
 
\usepackage{hyperref}
\hypersetup{colorlinks=false}
\usepackage{verbatim} 
\usepackage{xcolor}
\hypersetup{colorlinks=false,linkbordercolor=red,linkcolor=green,pdfborderstyle={/S/U/W 1}}

\usepackage[normalem]{ulem}
\usepackage{lipsum}
\usepackage{fancyhdr}
\usepackage{lastpage}

\usepackage[english]{babel}

\usepackage{bbm}
\usepackage{enumerate}
\usepackage[all]{xy}

\usepackage{hyperref}
\hypersetup{colorlinks=false}

\allowdisplaybreaks
\newtheorem{thm}{Theorem}[section]
\newtheorem{pro}[thm]{Proposition}
\newtheorem{lm}[thm]{Lemma}

\newtheorem{conj}[thm]{Conjecture}

\numberwithin{equation}{section}
\newtheorem{defn}[thm]{Definition}

\theoremstyle{remark}

\theoremstyle{remark}
\newtheorem{rem}[thm]{Remark}

\DeclareMathOperator*{\Irr}{Irr}

\DeclareMathOperator*{\Ind}{\textsf{Ind}}

\DeclareMathOperator*{\Res}{Res}

\DeclareMathOperator*{\Hom}{Hom}

\newcommand{\Gr}{{\rm Gr}}

\newcommand{\s}{\simeq}

\newcommand{\EE}{\mathbb{E}}

\newcommand{\CC}{\mathbb{C}}
\newcommand{\NN}{\mathbb{N}}

\newcommand{\G}{\mathbf{G}}

\newcommand{\bH}{\mathbf H}
\newcommand{\norm}{\mathcal{N}}

\def\bG{\bold G}

\DeclareMathOperator*{\SL}{SL}
\DeclareMathOperator*{\GL}{GL}

\title[Distinguished Representations for $\SL(n,F)$]{Distinguished Representations for $\SL(n,F)$}

\author{Kwangho Choiy}
\address{School of Mathematical and Statistical Sciences, Southern Illinois University, Carbondale, IL 62901-4408} 
\email{kchoiy@siu.edu}

\author{Shiv Prakash Patel}
\address{Department of Mathematics, IIT Delhi, Hauz Khas, New Delhi - 110016} 
\email{shivprakashpatel@gmail.com}


\date{\today}  
\begin{document}

\begin{abstract}
Let $F$ be a finite field, and let $\mathbb{E}$ be either a quadratic field extension $E/F$ or the split algebra $F \oplus F$.
We study distinguished representations of $\SL_{2n}(F)$ by the subgroup $H_{\flat} := \SL_{2n}(F) \cap \GL_n(\EE)$, which is a variation of the work of Anandavardhanan and Prasad on distinguished representations of $\SL_{n}(\EE)$ by the subgroup $\SL_n(F)$.
This is in a similar framework of our earlier work of a $p$-adic non-split variation of Anandavardhanan-Prasad over finite fields.
We give a formula for the dimension of the complex vector space $\Hom_{H_\flat}(\pi_\flat, \mathbbm{1})$ in terms of certain characters of $F^{\times}$, where $\pi_{\flat}$ is an irreducible representation which is also distinguished by $H_{\flat}$. 
\end{abstract}

\maketitle

\section{Introduction}
Given a subgroup $H$ of a group $G,$
an irreducible representation $\pi$ of $G$ is called $(H, \chi)$-distinguished  if $\Hom_{H}(\pi, \chi) \neq 0$, where $\chi$ is a character of $H$. 
In particular, for $\chi = \mathbbm{1}$ the trivial character of $H$, if $\Hom_{H}(\pi, \mathbbm{1}) \neq 0$, then we simply say that $\pi$ is $H$-distinguished.

In our earlier work \cite{choiypatel2025},
we studied distinguished representations  of the group $\SL_n(D)$ by the subgroup $\SL_n^{*}(E)$, where $D$ is a quaternion division algebra over a non-archimedean local field $F$ of characteristic zero and $E/F$ is a quadratic field extension and $\SL_n^{*}(E) = \SL_n(D) \cap \GL_n(E)$. 
We obtained a multiplicity formula for distinguished representations.
This was a non-split variation of the work of Anandavardhanan and Prasad,  who first studied $\SL_2(F)$-distinguished representations of $\SL_2(E)$ \cite{ap03} and later followed by $\SL_n(F)$-distinguished representations of $\SL_n(E)$  \cite{ap18}.

In a similar framework as our previous study, the present paper explores an analogous variation over a finite field $F$ with $q$ elements:  a study of distinguished representations of $\SL_{2n}(F)$ by a certain subgroup.
More precisely, let $\EE$ be  either a quadratic field extension $E/F$ or $F \oplus F$, which is a degree 2 algebra over $F$.
We consider the following algebraic groups over $F$,
\[
\bG={\GL}_{2n}, ~~~~~~~\bH = {\Res}_{\EE/F} {\GL}_n,  ~~~~~~~ \bG_\flat= {\SL}_{2n}.
\]
Define $G :=\bG(F) = \GL_{2n}(F)$, $H := \bH(F) = \GL_{n}(\EE)$ and $G_{\flat} := \bG_{\flat}(F) = \SL_{2n}(F)$.
Note that $H$ is either $\GL_n(E)$ or $\GL_n(F) \times \GL_n(F)$ depending on whether $\EE$ is $E$ or $F \oplus F$, respectively.
Set 
\[
H_{\flat} := G_{\flat} \cap H = \{ g \in {\GL}_n(\EE) : {\mathcal{N}_{\EE/F} }(\det(g))=1 \},
\]
where $\norm_{\EE/F} : \EE^{\times} \rightarrow F^{\times}$ stands for the norm map $\norm_{E/F}$ if  $\EE= E$, and $\norm_{\EE/F}(x,y) = xy$ if $\EE=F \oplus F$. 
We then have the following diagram:
\begin{align}
\xymatrix{
G_\flat \ar@{-}[d]_{\cup} \ar@{-}[r]^{\subset} &G \ar@{-}[d]^{\cup}\\
H_\flat \ar@{-}[r]^{\subset}          &H
}
\end{align} 
The main aim of the paper is to establish a multiplicity formula for the $H_{\flat}$-distinguished representations $\pi_{\flat}$ of the group $G_{\flat}$ i.e. a formula for the dimension of the $\CC$-vector space $\Hom_{H_{\flat}}(\pi_{\flat}, \mathbbm{1})$.
We will use the $(H, \chi)$-distinguished representation of $G$ for any character $\chi$ of $H$ studied by Dipendra Prasad \cite{prasad20}.

To state the main results, we introduce the following definitions.
For an irreducible representation $\pi$ of $G={\GL}_{2n}(F),$ we define 
\[
\begin{array}{cl}
X_{\pi} :=& \{{\alpha} \in \widehat{F^\times}: \pi~\text{ is }~ (H,\alpha \circ \mathcal{N}_{\EE/F} \circ \det)-\text{distinguished}\}, \\
 Z_\pi := & \{\chi \in \widehat{F^\times} : \pi \s \pi \otimes \chi \}.
\end{array}
\]

Let $\pi_\flat$ be an irreducible representation of $G_{\flat}$ which is $H_{\flat}$-distinguished. 
Then there exists an  irreducible representation $\pi$ of $G$ which is $H$-distinguished and  
$\pi_\flat \subset {\Res}^{G}_{G_\flat} (\pi)$ (see Lemma \ref{key pro}).
For $\alpha \in \widehat{\EE^{\times}}$ (here, $\EE^{\times}$ equals either $E^{\times}$ or $F^{\times} \times F^{\times}$), we set $m_{\alpha} = \dim_{\CC} {\Hom}_{H}( \pi, \alpha)$. 
The main result (Theorem \ref{mainthem}) is then the following.
\begin{thm} \label{mainthm:intro}
Let $\pi_{\flat}$ be a $H_{\flat}$-distinguished representation of $G_{\flat}$ and $\pi$ be an irreducible representation of $G$ such that $\pi_{\flat} \subset \Res^{G}_{G_{\flat}}(\pi)$ and $\pi$ is $H$-distinguished.
We have
\[
{\dim}_{\CC} {\Hom}_{H_{\flat}} (\pi_\flat,\mathbbm{1}) = \frac{ \sum_{\alpha \in X_\pi} m_\alpha}{|Z_\pi|}.
\]
In particular, if
$\pi_{\flat}$ is a cuspidal representation of $G_{\flat}$, then
\[
{\dim}_{\CC} {\Hom}_{H_{\flat}} (\pi_\flat,\mathbbm{1}) = \frac{|X_\pi|}{|Z_\pi|}.
\]
\end{thm}
We follow an approach similar to that of our earlier work in \cite{choiypatel2025} utilizing Prasad's  multiplicity one result for cuspidal representations (Theorem \ref{mult one of Prasad}) and classical results on restrictions of irreducible representations from a finite group to its normal subgroups (see Section \ref{restriction on finite groups}). To summarize, we make use of the intermediate subgroup 
\[
H^+ := G_\flat \cdot H  
\] 
which plays an important role here again. 
It turns out that any irreducible $H$-distinguished representation $\pi^+$ of $H^+$ contained in ${\Res}^{G}_{G\flat} (\pi)$ is unique and appears with multiplicity one in ${\Res}^{G}_{G_\flat} (\pi)$ (Lemma \ref{lemma1}). 
Since the multiplicity of $\pi^+$ upon restriction to $H$ is also one (Proposition \ref{fact}), the cardinality of those inequivalent components in $\pi^+|_{G_\flat}$ is equal to that of those in $\pi|_{\G_\flat}$ (\eqref{equalityincardinality})
We first see that $\dim_{\CC}{\Hom}_{H_{\flat}}(\pi,\mathbbm{1})$ equals
\begin{equation} \label{into(4.9)}
\dim_{\CC} {\Hom}_{G_\flat} (\pi, \pi_{\flat}) ~\cdot~ 
\dim_{\CC} {\Hom}_{H_{\flat}}(\pi_\flat,\mathbbm{1}) \cdot~ \big|\{\text{irreducible }\tau_\flat \text{ of }G_\flat: \tau_\flat \subset \pi^+|_{G_\flat}\}/\sim \big|.
\end{equation}
The cardinality of the last term in \eqref{into(4.9)} above turns out to be the number of characters of $F^{\times}$ stabilizing $\pi^+$, which is same as the number of character of $F^{\times}$ that stabilize $\pi$ due to Proposition \ref{fact}, which is $|Z_{\pi}|.$ 
Finally, we rewrite \eqref{into(4.9)} and use $\dim_{\CC}{\Hom}_{H_{\flat}}(\pi,\mathbbm{1}) = \sum_{\alpha \in X_\pi} m_\alpha$ (see \eqref{xyzusage}) to obtain the first part of Theorem \ref{mainthm:intro}.
In the case where $\pi_\flat$ is cuspidal, $\pi$ is also cuspidal (cf. \cite[Section 3]{leh72}). 
Then, by Prasad's Theorem \ref{mult one of Prasad}, we have $m_\alpha=1$ for all $\alpha \in X_{\pi}$, which completes the proof.

We remark that our main theorem yields a simple formula in terms of $X_{\pi}$ and $Z_{\pi}$ for cuspidal representations, thanks to the multiplicity one results of Prasad, i.e. $\dim_{\CC} \Hom_{H}(\pi, \chi) \leq 1$ for any cuspidal representation $\pi$ of $G$ and a character $\chi$ of $H$ for both cases $H=\GL_n(E)$ and $H=\GL_n(F) \times \GL_n(F)$. 
It should be noted that the multiplicity one result is not available to us for any irreducible representation $\pi$.
In fact, multiplicity one is not true for $H=\GL_n(F) \times \GL_n(F)$.
However, for $H=\GL_n(E)$, we expect the following.
\begin{conj} \label{conjecture}
Let $E/F$ be a quadratic extension. 
Let $\pi$ be an irreducible representation of $\GL_{2n}(F)$ and $\chi$ be a character of $E^{\times}$.    
Then
\[
\dim_{\CC} {\Hom}_{\GL_n(E)}(\pi, \chi \circ \det) \leq 1.
\]
\end{conj}
For cuspidal representation $\pi$, Conjecture \ref{conjecture} is a theorem of Prasad \cite[Theorem~4.1]{prasad20}.
In Section \ref{non-cuspidal-example}, we verify this conjecture for some non-cuspidal irreducible representations of $\GL_4(F)$.

\section*{Acknowledgements}
We thank Dipendra Prasad for helpful conversations and comments.
The first author is supported by a gift from the Simons Foundation (\#840755).
The second author acknowledges the support of the SERB MATRICS grant (MTR/2022/000782).

\section{Preliminaries}
\subsection{Notation}
Let $p$ be a prime and $q=p^n$ for some $n\in \NN.$ 
Throughout the paper, $F$ denotes a finite field with $q$ elements and $E/F$ be the quadratic extension. 
We write $\Irr(G)$ for the set of isomorphism classes of irreducible complex representations of a group $G.$
\subsection{Restriction from a finite group to its normal subgroup}  \label{restriction on finite groups}
Let $J_\flat$ and  $J$ be finite groups such that $J_\flat$ is a normal subgroup of $J$ and the quotient $J/J_\flat$ is abelian. 
Our primary references are \cite[(20.7) \& (20.8)]{jamesliebeck01}. 
Let $\delta_\flat \in \Irr(J_\flat)$ be given. We then have $\delta\in \Irr(J)$ such that $\delta_\flat \subset {\Res}_{J_\flat}^{J}(\delta).$
We write both $\Pi_{\delta_\flat}(J_\flat)$ and $\Pi_{\delta}(J)$ for the set of all nonequivalent irreducible representations in ${\Res}_{J_\flat}^{J}(\delta).$ 
We note that if $\Pi_{\delta_{\flat 1}}$ and $\Pi_{\delta_{\flat 2}}$
are any two irreducible constituents, then we have $\Pi_{\delta_{\flat 1}}(J) = \Pi_{\delta_{\flat 2}}(J)$.
Now we have the following decomposition
\[
{\Res}_{J_\flat}^{J}(\delta) = \langle \delta\rangle_{J_\flat} ~ \bigoplus _{\theta_\flat \in \Pi_{\delta}(J_\flat)}  \theta_\flat.
\]
Here, the multiplicity, denoted by $\langle \delta\rangle_{J_\flat},$ is the same for all $\theta_\flat \in \Pi_{\delta}(J_\flat)$. 
We define
\begin{equation*} \label{X(delta)}
I(\delta):= \{ \eta \in (J/J_\flat)^\vee : \delta\s \delta\eta \},
\end{equation*}
where $(J/J_\flat)^\vee = \Hom(J/J_\flat, \CC^{\times})$ is the set of all characters of $J/J_\flat.$ These are viewed as one-dimensional representations of $J$ that are trivial on the subgroup $J_{\flat}.$
For any $\delta_\flat \in \Pi_{\delta}(J_\flat),$ we consider the stabilizer of $\delta_{\flat}$ in $J$   
\[
J_{\delta_\flat}:= \{h \in J : {^h}{\delta_\flat} \s \delta_\flat
\}.
\]
We note that $J_{\delta_\flat}$ is a normal subgroup of $J,$ with $Z({J}) \cdot J_\flat \subseteq J_{\delta_\flat} \subseteq J,$ and the quotient $J/J_{\delta_\flat}$ acts by conjugation on the set $\Pi_{\delta}(J_\flat)$ simply and transitively.
We then obtain the equality 
\begin{equation} \label{cardinality}
|I(\delta)| = |\Pi_{\delta}(J_\flat)| \cdot \langle \delta_\flat, \delta\rangle_{J_\flat}^2.
\end{equation}
The following proposition is due to \cite[Proposition 3.23]{leh72} which will be utilized in several places below.
\begin{pro}[Lehrer] \label{fact}
Let $J$ be $\GL_n(F)$ and $J_\flat$ be $\SL_n(F).$ Then, for any $\delta \in \Irr(\GL_n(F)),$ we have 
\[
\langle \delta\rangle_{J_\flat} = 1.
\]
\end{pro}

\subsection{Restriction of representations of $G$ to $H$}
Our approach to obtain the multiplicity formula for the restriction of representations $\pi_{\flat}$ of $G_{\flat}$ to the subgroup $H_{\flat}$ relies on understanding the restriction of representations of $G$ to $H$. 
We begin by recalling some known results for the restriction of representations $\pi$ of $\GL_{2n}(F)$ to the subgroup $\GL_{n}(\EE)$. 
First, we state a theorem of Henderson \cite[Theorem 3.1.1 and Theorem 3.2.1]{Henderson2003}.
\begin{thm}[Henderson] \label{mult one of Hendeson}
Let $G= \GL_{2n}(F)$ and $H=\GL_{n}(\EE)$. Then $(G,H)$ is a Gelfand pair. 
In other words, for an irreducible representation $\pi$ of $\G$ we have 
\[
\dim_{\CC} {\Hom}_{H} (\pi, \mathbbm{1}) \leq 1.
\]
\end{thm}
The above theorem of Henderson was generalized by Prasad for any character $\chi$ of the subgroup $H$, but only for cuspidal representations of $G$ \cite[Theorem 4.1 and Theorem 5.3]{prasad20}. 
\begin{thm}[Prasad] \label{mult one of Prasad}
Let $G= \GL_{2n}(F)$ and $H=\GL_{n}(\EE)$.
Let $\pi$  be an irreducible cuspidal representation of $G$ and $\chi$ be a character of $H$. Then
\[
\dim_{\CC} {\Hom}_{H}(\pi, \chi) \leq 1.
\]
\end{thm}

\section{Multiplicity formula} \label{mainresult}
Recall that $G={\GL}_{2n}(F),  G_\flat={\SL}_{2n}(F)$,
\[ 
H = {\GL}_n(\EE) = \left\{ \begin{array}{ll}
     {\GL}_{n}(E) & \text{ if } \EE=E \\
     {\GL}_{n}(F) \times {\GL}_{n}(F) & \text{ if } \EE=F \oplus F. 
\end{array} \right. 
\]
and
\[
H_{\flat} := G_{\flat} \cap H = \{ g \in {\GL}_n(\EE) : {\mathcal{N}_{\EE/F} }(\det(g))=1 \}. 
\]
Then we have the following diagram:
\begin{center} 
\begin{tikzcd}[column sep={3em}] \label{slnD-slnL}
G_{\flat}={\SL}_{2n}(F)
  \arrow[dash]{rr}{\subset} 
  \arrow[dash]{dd}[swap]{\cup} 
  \arrow[dash]{rd}[inner sep=1pt]{\subset} 
& & 
G={\GL}_{2n}(F) 
 \arrow[shorten >= 10pt, dash]{dd}[inner sep=1pt]{\cup}
\\
& 
H^+ := {\SL}_{2n}(F) \cdot H
        \arrow[dash]{ur}[swap,inner sep=1pt]{\subset} 
& &

\\
H_{\flat}
  \arrow[swap,shorten >= 20pt, dash]{rr}{\subset} 
          \arrow[dash]{ur}[swap,inner sep=1pt]{\subset}
& & 
H
  \arrow[dash]{lu}[swap,inner sep=1pt]{\supset} 
\\
\end{tikzcd}
\end{center}

In this section we study the $H_{\flat}$-distinguished representations of $\SL_{2n}(F)$ following the framework established above and we obtain a multiplicity formula.
It should be noted that the embedding $H_{\flat} \hookrightarrow \SL_{2n}(F)$ is unique up to conjugation by elements of $\SL_{2n}(F)$ as we have
\[
|{\SL}_{2n}(F) \backslash {\GL}_{2n}(F)/ H|=1.
\]
\begin{rem} \label{someisoinfactorgroups}
We make the following remarks, which are pertinent to our setting.
\begin{enumerate}
\item We note that from ${\SL}_{2n}(F) \cdot H \cdot Z({\GL}_{2n}(F)) =: H^+$ we can skip writing $Z({\GL}_{2n}(F))$ as it is $F^{\times}$ which is already contained in $H$.

\item We have $H/H_{\flat} =H/(H \cap \SL_{2n}(F))= H^+ / \SL_{2n}(F),$ which is isomorphic to $F^\times$.  
This follows from the commutative diagram:
\begin{center} 
\begin{tikzcd}[column sep={3em}] \label{1}
\GL_{2n}(F)
  \arrow{r}{\det} 
&   
F^{\times} \\
\GL_{n}(\EE)
  \arrow{r}{\det}
  \arrow{u}{\cup}
&  
\EE^{\times} 
\arrow{u}{\norm_{\EE/F}}
\\
\end{tikzcd}
\end{center}
\end{enumerate}
\end{rem}

We begin with the following lemma which is a variation of our $p$-adic version \cite[Propostion 5.2]{choiypatel2025}.
\begin{lm} \label{key pro}
Let $\pi_{\flat} \in \Irr(\SL_{2n}(F))$ be $H_{\flat}$-distinguished.
Then there exists an $H$-distinguished representation $\pi \in \Irr({\GL}_{2n}(F))$ such that
\[
\pi_{\flat} \subset {\Res}^{{\GL}_{2n}(F)}_{\SL_{2n}(F)} ~ ({\pi}).
\]
\end{lm}
\begin{proof}
The main idea of the proof is similar to that of \cite[Propostion 5.2]{choiypatel2025}. Let $\pi \in \Irr({\GL}_{2n}(F))$ be  such that 
\[
\pi_\flat \subset {\Res}^{{\GL}_{2n}(F)}_{\SL_{2n}(F)} (\pi).
\]
For any $\sigma \in \Irr(H),$ the vector space of linear functionals
\[
{\Hom}(\pi, \sigma)
\]
is equipped with a natural $H$-action given by
\[
(g\cdot \lambda)(v)= \sigma(g)(\lambda (\pi(g^{-1}) \cdot v)~~~ \text{ for all } v \in \pi, ~~ g\in H.
\]
We now consider the \( H_{\flat} \)-invariant subspace
\[
({\Hom}(\pi, \sigma))^{H_{\flat}} = {\Hom}_{H_{\flat}}(\pi, \sigma).
\]
By taking $\sigma= \mathbbm{1}$, the trivial representation of $H_{\flat}$, the space ${\Hom}_{H_{\flat}}(\pi, \mathbbm{1})$ carries a natural \( H \)-action, making it a representation of $H$ on which $H_{\flat}$ acts trivially. 
Therefore, as a representation of \( H \), the space ${\Hom}_{H_{\flat}}(\pi, \mathbbm{1})$ is a direct sum of the characters of $H/H_{\flat} \cong F^{\times}$.
Thus, we can write 
\begin{equation} \label{m_alpha}
{\Hom}_{H_{\flat}}(\pi, \mathbbm{1}) \s \bigoplus_{\alpha \in \widehat{H/H_{\flat}}} m_{\alpha} \, \alpha,
\end{equation}
where $m_\alpha = \dim _{\CC}\Hom_H(\pi, \alpha).$ 

Since $0 \neq {\Hom}_{H_{\flat}}(\pi_{\flat}, \mathbbm{1}) \subset {\Hom}_{H_{\flat}}(\pi,\mathbbm{1})$, there exists ${\alpha} \in \widehat{H/H_{\flat}} \s \widehat{F^\times}$ (see Remark \ref{someisoinfactorgroups}(2)) such that $m_{{\alpha}} \geq 1$.
On the other hand, we have $F^\times \s \GL_{2n}(F)/\SL_{2n}(F) \cong F^{\times}$. 
Therefore, the irreducible representation $\pi \otimes {\alpha}^{-1}$ of $\GL_{2n}(F)$ is a $H$-distinguished and satisfies $\pi_{\flat} \subset {\Res}^{{\GL}_{2n}(F)}_{\SL_{2n}(F)} ~ (\pi \otimes {\alpha}^{-1})$. 
\end{proof}

\begin{lm} \label{lemma1}
Let $\pi$ be an irreducible representation of $\GL_{2n}(F)$ that is $H$-distinguished. 
Write
\[
{\Res}^{{\GL}_{2n}(F)}_{H^+} (\pi) \s   \bigoplus_{j} m_j ~ \pi^+_{j},
\]
where $\pi_j^{+} \in \Irr(H^{+})$ and $m_j = \dim \Hom_{H^{+}} (\pi, \pi_{j}^{+})$.
Then we have \( m_j=1 \) for all \( j \), and there is exactly one $\pi^+_{j}$ which is $H$-distinguished.
\end{lm}
\begin{proof}
Proposition \ref{fact} implies that all \( \pi_{j}^{+} \) have the same multiplicity, and that \( m_j = 1 \) for all \( j \). 
Since $H \subset H^+$, we have
\[
{\Hom}_{H} (\pi, \mathbbm{1} ) \s   \bigoplus_{j} {\Hom}_{H} ( \pi^+_j, \mathbbm{1}) .
\]
By Theorem \ref{mult one of Hendeson},
$\dim {\Hom}_{H} (\pi, \mathbbm{1} ) \leq 1$.
As $\pi$ is assumed to be $H$-distinguished, it follows that exactly one of the $\pi_{j}^{+}$ is $H$-distinguished.
\end{proof}

\begin{lm} \label{lemma2}
Let $\pi^{+} \in \Irr(H^+)$ be $H$-distinguished, and let $\pi_{\flat}, \pi_{\flat}' \in \Irr(\SL_{2n}(F))$ be contained in $\Res^{H^+}_{\SL_{2n}(F)} (\pi^{+})$. 
Then
\[
{\Hom}_{H_{\flat}} (\pi_{\flat},\mathbbm{1}) \cong {\Hom}_{H_{\flat}} (\pi_{\flat}',\mathbbm{1}).
\]
\end{lm}
\begin{proof}
This follows from the existence of an element \( g \in H \) such that \( \pi_{\flat}' \cong \pi_{\flat}^{g} \).
\end{proof} 

\begin{defn} \label{def XYZ}
Let $\pi$ be an irreducible representation of $\GL_{2n}(F)$. 
Define:
\[
\begin{array}{cl}
X_{\pi} :=&\{{\alpha} \in \widehat{F^\times}: \pi~\text{ is }~ (H,\alpha \circ \mathcal{N}_{\EE/F} \circ \det)-\text{distinguished}\}, \\
 Z_\pi := & \{\chi \in \widehat{F^\times} : \pi \s \pi \otimes \chi \}.
\end{array}
\]
\end{defn}

\begin{rem}
Note that $X_{\pi}$ need not be a subgroup of $\widehat{F^{\times}}$, but $Z_{\pi}$ is a subgroup. If $\alpha \in X_{\pi}$ and $\chi \in Z_{\pi}$, then $\alpha \chi \in X_{\pi}$. 
It follows that $X_{\pi}$ is a union of cosets of $Z_{\pi}$. 
Also, the group $\widehat{F^\times}$ of characters in $X_\pi$ can be identified with $\widehat{H/H_{\flat}}$  (refer to Remark \ref{someisoinfactorgroups}(2)).
\end{rem}

The following is the main theorem which is a variation of \cite[Proposition 3.6]{ap18}.

\begin{thm} \label{mainthem}
Let $\pi_\flat \in \Irr(\SL_{2n}(F))$ be $H_{\flat}$-distinguished. Choose an $H$-distinguished irreducible representation $\pi$ of ${\GL}_{2n}(F)$ such that  
\[
\pi_\flat \subset {\Res}^{{\GL}_{2n}(F)}_{\SL_{2n}(F)} (\pi).
\]
\begin{enumerate}
\item Let $\alpha \in \widehat{\EE^{\times}}$ and $m_{\alpha} = \dim_{\CC} {\Hom}_{H}( \pi, \alpha)$. 
Then
\[
{\dim}_{\CC} {\Hom}_{H_{\flat}} (\pi_\flat,\mathbbm{1}) = \frac{ \sum_{\alpha \in X_\pi} m_\alpha}{|Z_\pi|}.
\]
\item If $\pi_{\flat}$ is a cuspidal representation of $\SL_{2n}(F)$, then
\[
{\dim}_{\CC} {\Hom}_{H_{\flat}} (\pi_\flat,\mathbbm{1}) = \frac{|X_\pi|}{|Z_\pi|}.
\]
\end{enumerate}
\end{thm}

\begin{proof}
From \eqref{m_alpha}, we recall the isomorphism
${\Hom}_{H_{\flat}}(\pi, \mathbbm{1}) \cong \bigoplus_{\alpha \in \widehat{H/H_{\flat}}} m_{\alpha} \, \alpha$ as representation of $H/H_{\flat}$.  
Using Definition \ref{def XYZ}, this gives
\begin{equation} \label{xyzusage}
\dim_{\CC}{\Hom}_{H_{\flat}}(\pi,\mathbbm{1}) = \sum_{\alpha \in X_\pi} m_\alpha.
\end{equation}
Note that in the multiplicity-one cases of Theorem \ref{mult one of Prasad}, the right-hand side equals \( |X_\pi| \).

Given an $H_{\flat}$-distinguished $\pi_\flat \in \Irr(\SL_{2n}(F)),$ we choose an $H$-distinguished irreducible representation $\pi$ of ${\GL}_{2n}(F)$ such that  
\[
\pi_\flat \subset {\Res}^{{\GL}_{2n}(F)}_{\SL_{2n}(F)} (\pi).
\]
Using Lemma \ref{lemma1}, let $\pi^+ \subset \Res_{H^{+}}^{\GL_{2n}(F)} (\pi)$ be the unique $H$-distinguished representation of $H^+$ with multiplicity one.
Then we have the decomposition
\[
\pi^+|_{\SL_{2n}(F)} \s \langle \pi^+  \rangle_{\SL_{2n}(F)} 
\left( \bigoplus_{\{\tau_\flat \in \Irr({\SL}_{2n}(F)): \tau_\flat \subset \pi^+|_{\SL_{2n}(F)}\}/\sim}~\tau_\flat \right),
\]
where  $\langle \pi^+  \rangle_{\SL_{2n}(F)}$ denotes the multiplicity in the restriction (see Section \ref{restriction on finite groups}).
By Proposition \ref{fact}, we get 
\[
\langle \pi^+  \rangle_{\SL_{2n}(F)}=1.
\]
Hence, for any $\tau_{\flat} \in \Irr(\SL_{2n}(F))$ with $\tau_{\flat} \subset \Res_{\SL_{2n}(F)}^{\GL_{2n}(F)} (\pi)$, there exists a unique $\pi^{+} \in \Irr(H^{+})$ such that $\tau_{\flat} \subset \Res_{H^{+}}^{\GL_{2n}(F)} (\pi^{+})$.
Moreover, if $\tau_{\flat}$ is $H_{\flat}$-distinguished, then $\pi^{+}$ is $H$-distinguished.
We then have
\begin{equation} \label{equalityincardinality}
\big|\{\tau_\flat \in \Irr({\SL}_{2n}(F)): \tau_\flat \subset \pi^+|_{\SL_{2n}(F)}\}/\sim \big| = \big|\{\tau_\flat \in \Irr({\SL}_{2n}(F)): \tau_\flat \subset \pi|_{\SL_{2n}(F)}\}/\sim \big|.
\end{equation}
By Lemma \ref{lemma2}, we have
\begin{align}
\dim_{\CC}{\Hom}_{H_{\flat}}(\pi,\mathbbm{1}) 
=& \dim_{\CC}{\Hom}_{H_{\flat}}(\pi^{+},\mathbbm{1})   \nonumber \\
=&  \dim_{\CC} {\Hom}_{\SL_{2n}(F)} (\pi^+, \pi_{\flat}) ~\cdot~ 
\dim_{\CC} {\Hom}_{H_{\flat}}(\pi_\flat,\mathbbm{1})  \nonumber \\
 & ~ \cdot~ \big|\{\tau_\flat \in \Irr({\SL}_{2n}(F)): \tau_\flat \subset \pi^+|_{\SL_{2n}(F)}\}/\sim \big|.   \nonumber
\end{align}
Since $\dim_{\CC} {\Hom}_{\SL_{2n}(F)} (\pi^+, \pi_{\flat}) = 1$ by Proposition \ref{fact}, we obtain
\begin{equation}  \label{imp equation}
\dim_{\CC}{\Hom}_{H_{\flat}}(\pi,\mathbbm{1}) =   
\dim_{\CC} {\Hom}_{H_{\flat}}(\pi_\flat,\mathbbm{1})  \cdot~ \big|\{\tau_\flat \in \Irr({\SL}_{2n}(F)): \tau_\flat \subset \pi^+|_{\SL_{2n}(F)}\}/\sim \big|. 
\end{equation}
On the other hand, from \eqref{cardinality} and Proposition \ref{fact}, we have 
\begin{equation}  
\big|\{\tau_\flat \in \Irr({\SL}_{2n}(F)): \tau_\flat \subset \pi|_{\SL_{2n}(F)}\}/\sim \big|  = |\{\chi \in \widehat{({\GL}_{2n}/{\SL}_{2n}(F))} : \pi \otimes \chi \s \pi \}|. \nonumber
\end{equation}
By the group structure in Remark \ref{someisoinfactorgroups}(2) and the equality in \eqref{equalityincardinality}, we have
\begin{align} 
\big|\{\tau_\flat \in \Irr({\SL}_{2n}(F)): \tau_\flat \subset \pi^+|_{\SL_{2n}(F)}\}/\sim \big| &=
|\{\chi \in \widehat{F^\times} : \pi \otimes \chi \s \pi\}| \nonumber \\
&= |Z_{\pi}|.  \nonumber
\end{align}
Rewriting \eqref{imp equation} and using \eqref{xyzusage}, we find
\begin{align}
\dim {\Hom}_{H_{\flat}}(\pi,\mathbbm{1}) &= |Z_{\pi}| \cdot \dim {\Hom}_{H_{\flat}}(\pi_{\flat},\mathbbm{1}) \nonumber \\
\implies  \dim {\Hom}_{H_{\flat}}(\pi_{\flat},\mathbbm{1}) &= \dfrac{\dim {\Hom}_{H_{\flat}}(\pi,\mathbbm{1})}{|Z_{\pi}|} \nonumber \\
&= \dfrac{\sum_{\alpha \in X_{\pi}} m_{\alpha}}{|Z_{\pi}|} \nonumber
\end{align}
This proves part (1) of the theorem.
For part (2), as $\pi$ is cuspidal (c.f., \cite[Section 3]{leh72}), it follows from Theorem \ref{mult one of Prasad}, where we have $m_\alpha=1$ for $\alpha \in X_{\pi}.$ 
This completes the proof.
\end{proof}

\section{On multiplicity one for $\GL_4(F)$} \label{non-cuspidal-example}
Our multiplicity formula is in terms of $\dim {\Hom}_{\GL_n(E)} (\pi, \chi)$, which is less than or equal to one in the $p$-adic case \cite{lu23} and for cuspidal representations $\pi$ in the finite field case \cite{prasad20}. 
As mentioned in the introduction that we expect that the multiplicity one result remains valid for all the irreducible representations which we verify, in particular, for some irreducible representations of $\GL_4(F)$.

\subsection{Principal series representation induced from $(2,2)$ parabolic}
Consider the $P_{2,2}$ parabolic subgroup of $\GL_4(F)$, whose Levi subgroup $\GL_2(F) \times \GL_2(F)$ of $P_{2,2}$ is embedded in $\GL_4(F)$ as block diagonal matrices and the unipotent radical of $P_{2,2}$ is isomorphic to $M_2(F)$.
Let $\sigma_1, \sigma_2 \in \Irr(\GL_2(F))$ be irreducible representations and $\pi = \Ind_{P_{2,2}}^{\GL_4(F)}(\sigma_1 \otimes \sigma_2)$ a parabolically induced representation of $\GL_4(F)$.
We use Mackey theory to analyze the restriction $\Res^{\GL_4(F)}_{\GL_2(E)}(\pi)$, which requires a description of  the following double coset space
\[
{\GL}_2(E) \backslash {\GL}_4(F)/P_{2,2}.
\]
Note that $\GL_4(F)/P_{2,2}$ can be identified with the Grassmannian $\Gr(4,2)$, i.e. the set consisting of all 2-dimensional subspaces of a 4-dimensional vector space over $F$.
Note that $\GL_4(F)$ acts naturally on $\Gr(4,2)$ and hence the subgroup $\GL_2(E)$. 
Then $\Gr(4,2)$ is the union of disjoint $\GL_2(E)$-orbits and each orbit corresponds to a double coset in ${\GL}_2(E) \backslash {\GL}_4(F)/P_{2,2}$. 
We then have the following decomposition.

\begin{lm} \label{orbits for (2,2) parabolic}
Let $w \in \GL_4(F)$ be the permutation matrix obtained by interchanging the second and third rows of the identity matrix.
Then
\[
{\GL}_4(F) = \left( {\GL}_2(E) \cdot I \cdot P_{2,2} \right) \quad \bigsqcup \quad \left( {\GL}_2(E) \cdot w \cdot P_{2,2} \right).
\]
\end{lm}

\begin{proof}
Let us consider \( E \oplus E \) as a 4-dimensional vector space over \( F \). Define two 2-dimensional subspaces: 
\begin{itemize} 
\item \( W_1 = E \), embedded in $E \oplus E$ via \( x \mapsto (x, 0) \), 
\item \( W_2 = F \oplus F \), embedded in $E \oplus E$ via \( (x, y) \mapsto (x, y) \), using the inclusion \( F \subset E \). 
\end{itemize}

It can be easily observed that the stabilizer of $W_1$ is $B_2(E)$ the Borel subgroup of $\GL_2(E)$ consisting of upper triangular matrices and the stabilizer of $W_2$ is $\GL_2(F) \subset \GL_2(E)$. 
Evidently,  the stabilizer of $W_1$ is not isomorphic to the stabilizer of $W_2$.
Hence $W_1, W_2$ represent two distinct orbits, say $\mathcal{O}_{W_1}, \mathcal{O}_{W_2}$ in $\Gr(4,2)$ under $\GL_2(E)$-action.
It can be easily seen that $| \Gr(4,2)| = |\mathcal{O}_{W_1}| + |\mathcal{O}_{W_2}|$.
Then 
\[
\Gr(4,2) = \mathcal{O}_{W_1} \sqcup \mathcal{O}_{W_2}.
\]
and then the lemma follows by writing the representatives of the corresponding double cosets.
\end{proof}
By Mackey theory together with Lemma \ref{orbits for (2,2) parabolic}, we have
\begin{align}
\begin{array}{ll}
{\Res}_{\GL_2(E)}^{\GL_4(F)} \left( {\Ind}_{P_{2,2}}^{\GL_4(F)} (\sigma_1 \otimes \sigma_2) \right)
& \cong {\Ind}_{P_{2,2} \cap {\GL}_2(E)}^{\GL_2(E)}  (\sigma_1 \otimes \sigma_2) 
\bigoplus 
{\Ind}_{P_{2,2}^{w} \cap \GL_2(E)}^{\GL_2(E)} (\sigma_1 \otimes \sigma_2)^{w} \\
& \cong {\Ind}_{B_2(E)}^{\GL_2(E)}  (\sigma_1 \otimes \sigma_2) 
\bigoplus 
{\Ind}_{\GL_2(F)}^{\GL_2(E)} (\sigma_1 \otimes \sigma_2).
\end{array}
\end{align}
Therefore,
\begin{equation} \label{greatexamplegl4}
\begin{array}{ll}
{\Hom}_{\GL_2(E)}(\pi, \chi) 
& \cong {\Hom}_{\GL_2(E)} \left( {\Ind}_{B_2(E)}^{\GL_2(E)}  (\sigma_1 \otimes \sigma_2), \chi \right)  
\bigoplus
{\Hom}_{\GL_2(E)} \left( {\Ind}_{\GL_2(F)}^{\GL_2(E)} (\sigma_1 \otimes \sigma_2), \chi \right) \\
& \cong {\Hom}_{T(E)} \left( (\sigma_1 \otimes \sigma_2), \chi \right) \bigoplus
{\Hom}_{\GL_2(F)} \left( (\sigma_1 \otimes \sigma_2), \chi \right) \\
& \cong  {\Hom}_{E^{\times}} (\sigma_1, \chi) \otimes {\Hom}_{E^{\times}}(\sigma_2, \chi)
\bigoplus  
{\Hom}_{\GL_2(F)} \left( (\sigma_1 \otimes \sigma_2), \chi \right).
\end{array}
\end{equation}
The following is a simple observation:
\begin{equation} \label{sigma res to E}
{\Hom}_{E^{\times}}(\sigma_i, \chi) \neq 0 
\implies \chi|_{F^{\times}} = \omega_{\sigma_i}, \text{ the central character of } \sigma_i.
\end{equation}
Note that
\begin{align} \label{hom=1}
{\Hom}_{\GL_2(F)}(\sigma_1 \otimes \sigma_2, \chi \circ \det) \neq 0 \iff \sigma_2 \cong \sigma_1^{\vee} \otimes \chi \circ \det.
\end{align}
Now we have the following proposition.
\begin{pro}
Let $\sigma_1, \sigma_2$ be irreducible representations of $\GL_2(F)$ and, let $\pi = \Ind_{P_{2,2}}^{\GL_4(F)}(\sigma_1 \otimes \sigma_2)$.
If $\pi$ is irreducible then $\dim {\Hom}_{\GL_2(E)}(\pi, \chi) \leq 1$. 
\end{pro}
\begin{proof}
Suppose $\dim {\Hom}_{\GL_2(E)}(\pi, \chi) > 1$. 
We recall the following facts.
\begin{itemize}
\item From the character table of $\GL_2(F)$ (\cite[p.70]{Fulton-Harris}) it follows that
the pair $(\GL_2(F), E^{\times})$ is a strong Gelfand pair, i.e. $\dim {\Hom}_{E^{\times}}(\sigma, \chi) \leq 1$ for any irreducible representation $\sigma$ of $\GL_2(F)$.
\item From \cite[Lemma 2.4(c)]{gross1991}, it follows that the pair $(\GL_2(F) \times \GL_2(F), \triangle \GL_2(F))$ is also a strong Gelfand pair, i.e. $\dim {\Hom}_{\GL_2(F)}(\sigma_1 \otimes \sigma_2, \chi) \leq 1$ for any irreducible representations $\sigma_1, \sigma_2$ of $\GL_2(F)$.

\end{itemize}
Then by \eqref{greatexamplegl4}, if $\dim {\Hom}_{\GL_2(E)}(\pi, \chi) > 1$, then we must have
\[
{\Hom}_{E^{\times}} (\sigma_1, \chi) \cong \CC, ~~  {\Hom}_{E^{\times}}(\sigma_2, \chi) \cong \CC, ~~
{\Hom}_{\GL_2(F)} \left( (\sigma_1 \otimes \sigma_2), \chi \right) \cong \CC.
\]
By \eqref{sigma res to E}, for $i=1,2$ we have ${\Hom}_{E^{\times}} (\sigma_i, \chi) \cong \CC \implies \chi|_{F^{\times}} = \omega_{\sigma_i}$ the central character of $\sigma_i$. By \eqref{hom=1}, ${\Hom}_{\GL_2(F)} \left( (\sigma_1 \otimes \sigma_2), \chi \right) \cong \CC \implies \sigma_2 \cong \sigma_1^{\vee} \otimes \chi \circ \det \cong \sigma_1^{\vee} \otimes \omega_{\sigma_1} \circ \det \cong \sigma_1$. Since $\sigma_1 \cong \sigma_2$, by the irreducibility criterion of the parabolic induction,  $\pi = {\Ind}_{P_{2,2}}^{\GL_4(F)} (\sigma_1 \otimes \sigma_2)$ is reducible, contradicting the assumption.
\end{proof}
\subsection{Principal series representation induced from $(1,3)$ parabolic} 
Now we consider the other maximal parabolic subgroup of $\GL_4(F)$, namely the $(1,3)$-parabolic given by
\[
P_{1,3} = \left\{ \left( \begin{matrix} a & X \\ 0 & D \end{matrix} \right) : a \in F^{\times}, D \in {\GL}_3(F), X \in F^3 \right\}.
\]
Let $\mu : F^{\times} \to \CC^{\times}$ be a character and $\tau$ an irreducible cuspidal representation of $\GL_3(F)$.
Consider the principal series representation $\pi := {\Ind}_{P_{1,3}}^{\GL_4(F)}(\mu \otimes \tau)$. In order to understand $\pi|_{\GL_2(E)}$ using Mackey theory, we need a description of the double coset space
\[
{\GL}_2(E) \backslash {\GL}_4(F) / P_{1,3}.
\]
Note that $\GL_4(F) / P_{1,3}$ can be identified with the Grassmannian $\Gr(4,1)$, i.e. the set consisting of all 1-dimensional subspaces of a 4-dimensional vector space over $F$.
Then $\GL_2(E)$ acts naturally on $\Gr(4,1)$ via the natural action of $\GL_4(F)$. 
The following is a simple observation describing the above mentioned double coset space for which we omit the proof.
\begin{lm} \label{double coset for P13}
The action of $\GL_2(E)$ on $\Gr(4,1)$ is transitive. In other words, $\GL_4(F) = \GL_2(E) \cdot P_{1,3}$.    
\end{lm}
Using Mackey theory together with Lemma \ref{double coset for P13}, we get
\begin{align}
{\Res}^{\GL_4(F)}_{\GL_2(E)} (\pi) \cong {\Ind}_{\GL_2(E) \cap P_{1,3}}^{\GL_2(E)} (\mu \otimes \tau).  \nonumber
\end{align}
Then for any character $\chi : E^{\times} \to \CC^{\times}$, we get 
\begin{align}
{\Hom}_{\GL_2(E)}(\pi, \chi \circ \det)
& \cong {\Hom}_{\GL_2(E)}\left( {\Res}_{\GL_2(E)}^{\GL_4(F)}(\pi), \chi \circ \det \right) \nonumber \\
& \cong  {\Hom}_{\GL_2(E)} \left( {\Ind}_{H}^{\GL_2(E)} (\mu \otimes \tau), \chi \circ \det \right) \nonumber \\
& \cong {\Hom}_{H} \left( (\mu \otimes \tau)|_{H}, \chi \circ \det \right).  \nonumber
\end{align}
Let us describe $H := \GL_2(E) \cap P_{1,3}$. 
We fix an embedding $E \hookrightarrow M_2(F)$. 
Then, 
\begin{align}H 
= {\GL}_2(E) \cap P_{1,3} 
= \left\{ \left( \begin{matrix}  
a & 0 & b_1 & b_2 \\
0 & a & b_3 & b_4 \\
0 & 0 & d_1 & d_2 \\ 
0 & 0 & d_3 & d_4 
\end{matrix} \right) \in {\GL}_4(F) : a \in F^{\times}, 
\left( \begin{matrix}
d_1 & d_2 \\ d_3 & d_4
\end{matrix} \right) \in E^{\times},
\left( \begin{matrix}
b_1 & b_2 \\ b_3 & b_4
\end{matrix} \right) \in E
\right\}  \nonumber
\end{align}
Let $N \subset H$ be defined by
\begin{align}
N := \left\{ \left( \begin{matrix}  
1 & 0 & b_1 & b_2 \\
0 & 1 & b_3 & b_4 \\
0 & 0 & 1 & 0 \\ 
0 & 0 & 0 & 1 
\end{matrix} \right) : \left( \begin{matrix}
b_1 & b_2 \\ b_3 & b_4
\end{matrix} \right) \in E  
\right\}   \nonumber
\end{align}
By definition, it is clear that $\chi \circ \det$ is trivial on $N$.
On the other hand, 
\begin{align}
(\mu \otimes \tau) \left( \begin{matrix}  
1 & 0 & b_1 & b_2 \\
0 & 1 & b_3 & b_4 \\
0 & 0 & 1 & 0 \\ 
0 & 0 & 0 & 1 
\end{matrix} \right) 
= \tau \left( \begin{matrix}  1 & b_3 & b_4 \\  0 & 1 & 0 \\ 0 & 0 & 1 \end{matrix} \right)   \nonumber
\end{align}
Since $\tau$ is a cuspidal representation $\mu \otimes \tau$ restricted to $N$ has no trivial character, therefore 
\begin{align}
{\Hom}_{N} (\mu \otimes \tau, \chi \circ \det) = 0     \nonumber
\end{align}
and hence 
\begin{align}
{\Hom}_{H} (\mu \otimes \tau, \chi \circ \det) = 0.      \nonumber
\end{align}
Thus we have proved the following.
\begin{pro}
Let $\pi = {\Ind}_{P_{1,3}}^{\GL_4(F)}(\mu \otimes \tau)$, where $\mu$ is a character of $F^{\times}$ and $\tau$ is an irreducible cuspidal representation of $\GL_3(F)$. Then $ {\Hom}_{\GL_2(E)}(\pi, \chi) =0$.    
\end{pro}


\end{document}